\documentclass[12pt]{article} 
\usepackage[utf8]{inputenc}
\usepackage[T2A]{fontenc}\usepackage[english]{babel}
\usepackage{amsthm, amsmath, amssymb, amsbsy, amscd, amsfonts, latexsym, euscript,epsf,stackrel}
\usepackage{authblk}
\usepackage{graphicx} 
\usepackage{epsfig}

\usepackage{tikz}
\usetikzlibrary{positioning,arrows}

\newtheorem{theorem}{Theorem}[section]

\newtheorem{lemma}[theorem]{Lemma}
\newtheorem{proposition}[theorem]{Proposition}
\theoremstyle{definition}

\newtheorem{question}[theorem]{Question}

\newtheorem{definition}[theorem]{Definition}
\newtheorem{example}[theorem]{Example}
\newtheorem{remark}[theorem]{Remark}
\numberwithin{equation}{subsection}
\newtheorem*{ack}{Acknowledgement}
\usepackage[all,cmtip]{xy}
\usepackage{xcolor}
\usepackage{xfrac}
\usepackage{pb-diagram}
\usepackage{graphicx}
\usepackage{float}

\newcommand{\Inn}{\operatorname{Inn}}

\newcommand{\T}{\operatorname{T}}

\newcommand{\id}{\mathrm{id}}

\makeatletter
\newcommand{\subjclass}[2][1991]{%
  \let\@oldtitle\@title%
  \gdef\@title{\@oldtitle\footnotetext{#1 \emph{Mathematics subject classification.} #2}}%
}
\newcommand{\keywords}[1]{%
  \let\@@oldtitle\@title%
  \gdef\@title{\@@oldtitle\footnotetext{\emph{Key words.} #1.}}%
}
\makeatother

\def\bea{\begin{eqnarray}}
\def\eea{\end{eqnarray}}
\def\nn{\nonumber}

\begin{document}
\title{Self-distributive algebras and bialgebras}

\author[1,2,3]{Valeriy G. Bardakov\footnote{bardakov@math.nsc.ru}}
\author[3]{Tatiana A. Kozlovskaya\footnote{konus\_magadan@mail.ru}}
\author[4,5]{Dmitry~V.~Talalaev\footnote{dtalalaev@yandex.ru}}
\affil[1]{\small Sobolev Institute of Mathematics, Novosibirsk 630090, Russia.
}
\affil[2]{\small Novosibirsk State Agrarian University, Dobrolyubova street, 160, Novosibirsk, 630039, Russia.
}
\affil[3]{\small Regional Scientific and Educational Mathematical Center of Tomsk State University,
36 Lenin Ave., 14, 634050, Tomsk, Russia.  
}
\affil[4]{\small Lomonosov Moscow State University, 
119991, Moscow, Russia.
}
\affil[5]{\small Center of Integrable Systems, Demidov Yaroslavl State University, Yaroslavl, Russia, 150003, Sovetskaya Str. 14
}

\date{\today}

%\subjclass[2010]{Primary 17D99; Secondary 57M27, 16S34, 20N02}

\maketitle

\begin{abstract}
This article is devoted to the study of self-distributive algebraic structures: algebras, bialgebras; additional structures on them, relations of these structures with Hopf algebras, Lie algebras, Leibnitz algebras etc. The basic example of such structures are rack- and quandle bialgebras. But we go further - to the general coassociative comultiplication. The principal motivation for this work is  the development of the linear algebra related with a notion of a quandle in analogy with the ubiquitous role of group algebras in the category of groups with perspective applications to the theory of knot invariants. We give description of self-distributive algebras and show that some quandle algebras and some Novikov algebras are self-distributive. Also, we give a  full classification of counital self-distributive bialgebras in dimension 2 over~$\mathbb{C}$.
\end{abstract}

\tableofcontents

\section{Introduction and Preliminaries}
\label{Intro}
\medskip

A set-theoretical  solution  of the  Yang--Baxter equation (YBE) \cite{Drinfeld} on a set $X$ is a map $$R \colon X \times X \to X\times X$$ satisfying 
\begin{equation}
\label{YB}
	R_{12} R_{23} R_{12} = R_{23} R_{12} R_{23},
\end{equation}
 where   $R_{ij}  \colon X \times X \times X  \to  X \times X \times X$ acts as $R$ on the $i$-th and $j$-th factors and as the identity on the remaining factor. 
The pair $(X, R)$ is said to be a solution of the YBE or simply a  solution.

The YBE appears  in the papers of C.~N.~Yang  \cite{Yang} and R.~J.~Baxter \cite{Bax}.   
It is   one of the basic  equations in mathematical physics and in  low dimension topology.
It lies in the foundation of the theory of quantum groups \cite{Kas}, solvable models of statistical mechanics,   knot theory, braid theory \cite{KT}. 
 Faddeev, Reshetikhin and Takhtajan \cite{FRT} provided a way to obtain some bialgebras by means of $R$--matrices, where $R$ is a matrix solutions for the Yang--Baxter equation. Jones \cite{J} proposed a ``baxterization'' technique which
shows a connection between representations of the braid group and solutions of the Yang--Baxter
equation.

Set-theoretical solutions of the
YBE are closely related with many algebraic systems. One of them is a rack that is a non-empty set with one binary algebraic operation with  
 two axioms which provide an algebraic interpretation of Reidemeister moves. This structure also plays an important
role in the study of the YBE. 
Racks  are responsible for complete sequence of computable invariants for framed links and for 3-manifolds (see~\cite{FR}).   Andruskiewitsch and M. Gra\~na \cite{AG} proved that there is a one-to-one corresponding between racks and certain set-theoretic
solutions of the YBE.  
Racks are closely related to Leibniz algebras. Kinyon \cite{Kin} showed that a Leibniz algebra defines  a rack.

The construction of a groupoid algebra is well known. If $X$ is a groupoid, that is a non-empty set with one binary algebraic operation and $\Bbbk$ is a commutative associative ring with unit, then the set of finite formal linear combinations form an  algebra $\Bbbk[X]$. If $(X, \cdot)$ is an associative groupoid, then the algebra $\Bbbk[X]$ is associative.
This follows from the fact that the  associativity identity $(ab)c = a(bc)$ is tri-linear. 
Now, suppose that $(X, *)$ is a right self-distributive groupoid that  means that operation $*$ is self-distributive, 
$$
(u*v)*w = (u*w)*(v*w),
$$ 
for any $u, v, w \in X.$ Even if there is a linear structure on $X$ this condition is not tri-linear. 
%We will show that, for example,  the quandle algebra $\Bbbk[Q]$ is not self-distributive under the operation which is induced from the operation in~$Q$.
In the present article we are studying  self-distributive algebras and bialgebra, generalizing this identity. 
In \cite{CCES} there were defined self-distributive structures in the categories of coalgebras and cocommutative
coalgebras. The first examples of these structures are given by  vector spaces based on the elements of finite
quandles, the direct sum of a Lie algebra with its ground field, a Hopf algebra. The self-distributive operations of these structures provide solutions for the Yang--Baxter equation, and,
conversely, Yang--Baxter maps can be used to construct self-distributive
operations in certain categories.

In \cite{ABRW} there was defined a certain self-distributive multiplication on coalgebras, which leads to so-called rack bialgebra. In particular, it was showed that a rack algebra $\Bbbk[X]$ is a self-distributive bialgebras for the group-like comultiplication $\Delta(x) = x \otimes x$, $x \in X$.

Basic Lie theory gives the fundamental links between associative algebras, Lie algebras and groups. Some of these links are the passage
from an associative algebra $A$ to its underlying Lie algebra $A^{Lie}$ which is the vector space $A$ with the bracket $[a, b] = ab - ba$. On the other hand, to any
Lie algebra $\mathfrak{g}$ one may associate its universal enveloping algebra $U(\mathfrak{g})$ which is associative. Groups arise as groups of units in associative algebras. To any group $G$, one may associate its group algebra $\Bbbk[G]$ which is associative.

In an attempt to bring ring theoretic techniques to the study of quandles, a theory
of quandle rings analogous to the classical theory of group rings has been proposed
in \cite{BPS}, where several relations between quandles and their associated quandle
rings have been explored. 
Quandle rings of non-trivial quandles are non-associative, and it has
been proved in \cite{EFT}  that these rings are not even power-associative, which is the other
end of the spectrum of associativity.

Units in group rings play a fundamental role in the structure theory of group rings. Since each element of a quandle
 is an idempotent of the quandle ring,  idempotents play the same role as invertible  elements  in group rings.
 In \cite{ENS}
idempotents of quandle rings have been used for constructing proper enhancements of
the well-known quandle coloring invariant of knots and links in the 3-space.
Idempotents of integral quandle
rings of all three element quandles and the dihedral quandle of order four have been
computed in \cite{BPS-1}. 
 In  \cite{ENSS} there were  investigated idempotents in quandle rings, in particular, idempotents in quandle rings of
free products including free quandles.

The paper is organized as follows. In Section \ref{Intro} we recall some definitions of coalgebras, Hopf algebras, racks and quandles, which can be found in \cite{Kas, Joyce, Mat, CP}. 

In Sections \ref{self-dis} we introduce self-distributive algebras, self-distributive bialgebras, analize relations with other construction. In Section \ref{self-rel} we realize multiplie relations of self-distributive structures with Yang-Baxter equation, Novikov algebras, generalized  Jordan algebras, Lie algebras and Hopf algebras.. 

In Sections \ref{class} we give a  full classification of counital self-distributive bialgebras in dimension 2 over the field $\mathbb{C}$. This result extends \cite[Section 3.4]{CCES} where the partial result was obtained.

\bigskip

\subsection{Coalgebra, bialgebra,  and Hopf algebra} \label{coalg}
Let $C$ be a module over a commutative associative unital ring $\Bbbk$. A linear map 
$$
\Delta \colon C \to C \otimes_\Bbbk C = C \otimes C
$$
is called a comultiplication on $C$. Let us denote a comultiplication  of $a \in C$ by
\bea
\Delta(a)=\sum_i a_i^{(1)} \otimes a_i^{(2)} =a^{(1)}\otimes a^{(2)},\nn
\eea
where we use  Sweedler notations. 
A comultiplication is called coassociative  if
 $$
(\Delta \otimes \id_C )   \Delta  = (\id_C \otimes \Delta )   \Delta. 
 $$ 
In this case  the pair $(C, \Delta)$ is called a  coalgebra over $\Bbbk$. The coalgebra $(C, \Delta)$ is called cocommutative  if $\tau   \Delta =  \Delta$, where
$\tau \colon C \otimes C \to  C \otimes C$ is the canonical flip map $\tau(a \otimes b) = b  \otimes a$. A linear map $\varepsilon \colon C \to \Bbbk$  is called a counit for the coalgebra $(C, \Delta)$ if 
 $$
(\varepsilon \otimes \id_C)   \Delta  = (\id_C \otimes \varepsilon)   \Delta = \id_C, 
 $$ 
where we assume  $\Bbbk \otimes C \cong C$.

The triple $(C, \Delta, \varepsilon)$ is called a counital coalgebra. Moreover, a counital coalgebra $(C, \Delta, \varepsilon)$ equipped with an element ${\bf 1}$ is called coaugmented if $ \Delta({\bf 1}) = {\bf 1} \otimes {\bf 1}$ and $\varepsilon({\bf 1}) = 1 \in \Bbbk$.
In \cite{ABRW}  a coassociative, counital, coaugmented coalgebra is called a 
 $C^3$-coalgebra. In case the $C^3$-coalgebra is in addition cocommutative, it is called by a $C^4$-coalgebra.

A bialgebra $(H; \cdot, \Delta, \varepsilon)$  is a $\Bbbk$-module $H$,  $(H; \cdot)$ is an algebra (may be non-associative and without unit) and $(H;  \Delta, \varepsilon)$ is a coassociative  coalgebra with compatible  multiplication and comultiplication,  that means that the multiplication is a homomorphism of coalgebras and the comultiplication is a homomorphism of algebras. Both conditions are equivalent to the following
\bea
(hg)^{(1)}\otimes (hg)^{(2)}=h^{(1)} g^{(1)}\otimes h^{(2)} g^{(2)}.\nn
\eea 
An {\it antipode} of a bialgebra $H$ is an anti-automorphism of both the algebra and the coalgebra structure in H,
$$
S \colon H \to H
$$
such that
$$
m (S \otimes \id_H ) \Delta  = m (\id_H \otimes S) \Delta = \varepsilon \cdot 1_H.
$$

A Hopf algebra is an associative  bialgebra with a unit $1_H$,
a counit $\varepsilon$  and an antipode $S$.
 It also follows from the axioms of the Hopf algebra  that 
 $$
 S(1_H) = 1_H~\mbox{and}~\varepsilon \, S = \varepsilon \colon H \to \Bbbk.
 $$

\subsection{Racks and quandles}  \label{RQ}
A {\it quandle} is a non-empty set $Q$ with a binary operation $(x,y) \mapsto x * y$ satisfying the following axioms:
\begin{enumerate}
\item[(Q1)] $x*x=x$ for all $x \in Q$,
\item[(Q2)] For any $x,y \in Q$ there exists a unique $z \in Q$ such that $x=z*y$,
\item[(Q3)] $(x*y)*z=(x*z) * (y*z)$ for all $x,y,z \in Q$.
\end{enumerate}
~\\
An algebraic system satisfying only (Q2) and (Q3)  is called a {\it rack}. An algebraic system satisfying only (Q3)  is called a {\it self-distributive groupoid or shelf}.

If $Q$ is a rack, then from the second axiom it follows that if $c * a = b$, then one can define an operation $\bar{*} \colon Q \times Q \to Q$  by the rule $c = b \, \bar{*} \, a$. Hence, we can define the second rack axiom in the form:

(Q2') For any $a, b \in Q$ it holds
$$
(b \, \bar{*} \, a) * a = b =(b * a) \, \bar{*} \, a.
$$

A quandle  $Q$ is called {\it trivial} if $x*y=x$ for all $x, y \in Q$.  Unlike groups a trivial quandle can have arbitrary number of elements. We denote the $n$-element trivial quandle by $\T_n$ and an arbitrary trivial quandle by $\T$.
\medskip
~\\ 
Notice that the axioms (Q2) and (Q3) are equivalent to the map $S_x \colon Q \to Q$ given by $$S_x(y)=y*x$$ being an automorphism of $Q$ for each $x \in Q$. These automorphisms are called {\it inner automorphisms}, and the group generated by all such automorphisms is denoted by $\Inn(X)$.

\bigskip

\section{Self-distributive structures: definitions and examples} \label{self-dis}

\medskip 

%\subsection{Self-distributive algebras and self-distributive bialgebras} \label{SDA-SDC}
Here we recall the donor algebraic structure - the self-distributive bialgebra. But we start with the more elementary one.
An algebra $A$ over a field $\Bbbk$ is said to be self-distributive algebra if
$$
(a b) c = (a c) (b c),~~~a, b, c \in A.
$$

Recall \cite{CCES} (see  also \cite{ABRW}) that a counital  bialgebra $A$ with a comultiplication $\Delta$ is said to be a {\it self-distributive bialgebra} if  
$$
(a b) c = (a c^{(1)})  (b c^{(2)})
$$
for any $a, b, c \in A$

\medskip

\subsection{Rack coalgebras and  bialgebras }

Rack algebra and quandle algebra were introduced in \cite{BPS}.
Let $X$ be a rack and  $\Bbbk[X]$ the set of all formal finite $\Bbbk$-linear combinations of elements of $X$, that is,
$$\Bbbk[X]:=\Big\{ \sum_i\alpha_i x_i~|~\alpha_i \in \Bbbk,~ x_i \in X \Big\}.$$
Then $\Bbbk[X]$ is an additive abelian group with coefficient-wise addition. Define multiplication in $\Bbbk[X]$ by setting $$\big(\sum_i\alpha_i x_i\big) \big(\sum_j\beta_j x_j\big):=\sum_{i,j}\alpha_i\beta_j (x_i x_j).$$
Clearly, the multiplication is distributive with respect to addition from both left and right, and $\Bbbk[X]$ forms  a $\Bbbk$-algebra, which is called  the  {\it rack algebra} of $X$ with coefficients in the ring $\Bbbk$. Since $X$ is non-associative, unless it is a trivial rack, it follows that $\Bbbk[X]$ is a non-associative algebra in general. 
\par

Define the augmentation map
$$\varepsilon \colon \Bbbk[X] \to \Bbbk$$
 by setting $$\varepsilon \big(\sum_i\alpha_i x_i\big)= \sum_i\alpha_i .$$
Clearly, $\varepsilon$ is a surjective algebra homomorphism, and $I_\Bbbk(X):= \ker(\varepsilon)$ is a two-sided ideal of $\Bbbk[X]$, called the {\it augmentation ideal} of $\Bbbk[X]$. It is easy to see that $\{x-y~|~x, y \in X \}$ is  a generating set for $I_\Bbbk(X)$ as an $\Bbbk$-module. Further, if $x_0 \in X$ is a fixed element, then the set $\big\{x-x_0~|~x \in X \setminus \{ x_0\} \big\}$ is a basis for $I_\Bbbk(X)$ as a $\Bbbk$-module. 
\par
Let us define a comultiplication $\Delta \colon \Bbbk[X] \to \Bbbk[X] \otimes \Bbbk[X]$ by the rule $\Delta(x_i) = x_i \otimes x_i$ for any $x_i \in X$ and extend it to $\Bbbk[X]$ by linearity:
\bea
\Delta \big(\sum_i\alpha_i x_i\big)= \sum_i\alpha_i x_i\otimes x_i.\nn
\eea
Then $(\Bbbk[X], \Delta, \varepsilon)$ is a cocommutative coalgebra with the counit.

Further, take the extension $\Bbbk \oplus \Bbbk[X]$. Its elements
$$
\alpha + \sum_i\alpha_i x_i,~~\alpha, \alpha_i \in \Bbbk,~ x_i \in X. 
$$
Extend $\Delta$ and $\varepsilon$ to $\Bbbk \oplus \Bbbk[X]$ by linearly extending $\Delta(1) = 1 \otimes 1$ and $\varepsilon (1) = 1$, where $1 \in \Bbbk$. For arbitrary elements
\bea
\Delta \big(\alpha + \sum_i\alpha_i x_i\big) = \alpha (1 \otimes 1) + \sum_i\alpha_i x_i\otimes x_i,~~
\varepsilon \big(\alpha + \sum_i\alpha_i x_i\big) = \alpha  + \sum_i\alpha_i.\nn
\eea
We get a cocommutative coalgebra with the counit $(\Bbbk \oplus \Bbbk[X],  \Delta,  \varepsilon)$. 

Extend the multiplication from $\Bbbk[X]$ to $\Bbbk \oplus \Bbbk[X]$ by the rules
$$
1 \cdot x_i = 1,~~x_i \cdot 1 = 0,~~1 \cdot 1 = 0,
$$
we get an algebra $(\Bbbk \oplus \Bbbk[X], \cdot, \Delta,  \varepsilon)$ with multiplication and comultiplication. Remark that in general case the algebra $(\Bbbk \oplus \Bbbk[X], \cdot)$ is not self-distributive.

The following proposition can be found in  \cite[Proposition 3.1]{CCES}.
\begin{proposition}[\cite{CCES}]
$(\Bbbk \oplus \Bbbk[X], \cdot, \Delta,  \varepsilon)$ is a  counital self-distributive bialgebra.
\end{proposition}

\bigskip

\subsection{Quandle coalgebra and quandle  bialgebra} \label{QBA}

In this section we examine an algebra $\Bbbk[Q]$ if $Q$ is not only a rack but a quandle. 
	
\begin{lemma}
Suppose that a rack $(X, *)$ contains element $e$ such that $e * x = x * e = x$ for any $x \in X$. Then $X = \{e\}$ is a trivial 1-element rack.
\end{lemma}

\begin{proof}
From the definition of rack it follows that for any $a, b \in X$ the following holds
$$
(e * a) * b = (e * b) * (a * b).
$$
Since $e * a = a$, we have  $a * b =  b * (a * b)$. From the definition of the unit element it is equivalent to $e * (a * b) =  b * (a * b)$. The rack axiom  gives that $e = b$ for any $b \in X$.
\end{proof}

\medskip

If $(Q, \cdot)$ is a quandle, then  $(\Bbbk[Q], \cdot, \Delta,  \varepsilon)$ is a  counital self-distributive bialgebra.
Furthermore,  any element in $Q$ is an idempotent, but for 
$\Bbbk[Q]$ it is not true. In \cite{ABRW} it was defined a generalized idempotent.
\begin{proposition}
In a quandle bialgebra $(\Bbbk[Q], \Delta,  \varepsilon, *)$  any element $c$ is  a generalized idempotent that is
 $$
\sum c^{(1)}_i c^{(2)}_i = c.
$$
\end{proposition}

\begin{proof}	
Let 
$$
c = \alpha e + \sum_i \beta_i a_i,~~~\alpha, \beta_i \in \Bbbk
$$ 
be an element of $\Bbbk^\circ[Q]$. Then
$$
\Delta(c) = \alpha e \otimes e + \sum_i \beta_i a_i \otimes a_i,~~~\alpha, \beta_i \in \Bbbk,
$$
and
$$
\sum c^{(1)}_i c^{(2)}_i = \alpha e  e + \sum_i \beta_i a_i  a_i = c
$$
\end{proof}
	
Indeed, one could find a generalized invertibility axiom for rack and quandle bialgebras.
\begin{proposition}\label{invax}
Let $X$ be a rack, then we define a bilinear map 
\bea
\phi:\Bbbk[X]\otimes \Bbbk[X]\rightarrow \Bbbk[X]\nn
\eea
by acting on generators $a,b\in X$
\bea
\phi(a,b)=c\nn
\eea
where $c$ is a unique element of $X$ such that $c a=b.$ Such an element exists in virtue of the second rack axiom. Then
\bea
\label{eq_invax}
\phi(x^{(1)},y)x^{(2)}=y;\qquad \forall x,y\in\Bbbk[X].
\eea
\end{proposition}

\begin{proof}
This could be shown by direct calculation as in Proposition \ref{trivqv}. But we always could use the linearity argument, both sides of the equation \ref{eq_invax} are linear on arguments $x,y$, then \ref{eq_invax} could by verified just on basis elements, those are elements of the rack.
\end{proof}

\subsection{Self-distributive algebras}
As we know, the class of associative algebras is very large. For example, if $S$ is a semi-group, then the set $\Bbbk[S]$ of finite linear combinations  of elements of $S$ with coefficients in $\Bbbk$ form an associative algebra. But, if $Q$ is a rack, then the set $\Bbbk[Q]$ of finite linear combinations forms the rack algebra which is not self-distributive, it means that the following axiom
$$
(a b) c = (a c) (b c)
$$
does not holds for all $a, b, c \in \Bbbk[Q]$.
In the present section we are studying the  following question.
What can we say on  self--distributive algebras?

Using some ideas on linearization of identity (see \cite[Chapter 1]{ZSSS}), we prove that the class of self-distributive algebras is poor. 

Sasha Panasenko suggested to prove the following proposition.

\begin{proposition} \label{psd}
Let $A$ be a self--distributive algebra over a field $\Bbbk$ of characteristic $char(\Bbbk)\ne2$. Then $A$ satisfies the identity  $(AA) A = 0$.  
\end{proposition}

\begin{proof}

The self--distributivity means
$$
(a b) c = (a c) (b c)
$$
for any $a, b, c \in A$. If we put $c = d + f$ for some $d, f \in A$, then we get
$$
(a f) (b d) + (ad) (bf) = 0.
$$
If $d = f$, then $2 (ad) (bd) = 0$. Since $char(\Bbbk)\ne2$, then  $(ad) (bd) = 0$ and from self-distributivity follows $(a b) d = 0$ for any $a, b, d \in A$.
\end{proof}

The following statement shows that over the field of $char(\Bbbk)=2$ Proposition \ref{psd} does not hold.

\begin{proposition}
A quandle algebra $\mathbb{Z}_2[T]$ of a trivial quandle $(T; \cdot)$ is self-distributive.
\end{proposition}

\begin{proof}
We have to check the equality
$$
(a b) c = (a c) (b c)
$$
for any $a, b, c \in \mathbb{Z}_2[T]$. The left hand side of this  equality is
$$
(a b) c = (a \varepsilon(b)) c = a \varepsilon(b) \varepsilon(c).
$$
 The right hand side is
$$
(a c) (b c) = (a \varepsilon(c))  (b \varepsilon(c))  = a \varepsilon(b) (\varepsilon(c))^2.
$$
Hence, we need to check $\varepsilon(c) = (\varepsilon(c))^2$. This is always true for $\varepsilon(c)  \in \mathbb{Z}_2.$
\end{proof}

\medskip

As the other example of self-distributive algebra we give a Novikov algebra.

Recall that vector space $A$ with a bilinear operation $x,y\mapsto xy$ is called a Novikov algebra if the associator
\bea
(x,y,z)=(xy)z-x(yz)\nn
\eea
is {\it left-symmetric}:
\bea
(x,y,z)=(y,x,z)\nn
\eea
and the following property holds
\bea
(xy)z=(xz)y.\nn
\eea

\begin{example}
1) The algebra $(A_1 = \Bbbk e_1 \oplus \Bbbk e_2, \cdot)$ with multiplication
$$
e_1 \cdot e_1 = e_2 \cdot e_1 =e_2 \cdot e_2 = 0, ~~e_1 \cdot e_2 = e_2
$$
is a Novikov algebra. Since $A_1^2 A_1 = 0$, for any field $\Bbbk$ the  algebra $(A_1, \cdot)$ is a self-distributive algebra.

2) The algebra $(A_2 = \Bbbk e_1 \oplus \Bbbk e_2, \cdot, \Delta)$ with multiplication
$$
e_1 \cdot e_1 = e_1,~~e_2 \cdot e_1 =e_2,~~e_1 \cdot e_2 = e_2 \cdot e_2 = 0,
$$
is also a Novikov algebra, but $A_2^2 A_2 \not= 0$. On the other hand, it is not difficult to prove that 
for 2-element field $\Bbbk = \mathbb{Z}_2$ the  algebra $(A_2, \cdot)$ is a self-distributive algebra.
\end{example}

It is natural to formulate a question: what quandle axioms hold in algebras, in particular, in quandle algebras. In introduction we discussed the existence of idempotents in quandle algebras. In this subsection we discuss self-distributive algebras.
Now we proceed to the interpretation of the second axiom.

\begin{question}
Let $\Bbbk[X]$ be a rack algebra. For what elements $a, b\not= 0 \in \Bbbk[X]$ the equation 
\begin{equation} \label{eq}
x a = b,~~a, b\not= 0 \in \Bbbk[X]
\end{equation}
 has a solution in $\Bbbk[X]$?
\end{question}

For trivial quandles, the answer to this question gives

	\begin{proposition} \label{trivqv}
If $T$ is a trivial quandle, and 	$\Bbbk[T]$ is a quandle algebra. Then  an equation (\ref{eq}) has a solution if and only if $a$ does not lie in 
the augmentation ideal $I_\Bbbk(T)$.
	\end{proposition}

\begin{proof}
Let $T = \{ t_0, t_1, \ldots \}$. Then  the augmentation ideal is generated by elements 
$$
\tau_1 = t_1 - t_0, ~~\tau_2 = t_2 - t_0, \ldots.
$$
 If we take in  $\Bbbk[T]$ a new basis $t_0$, $\tau_1$, $\tau_2$, $\ldots$,
then the multiplication is defined by the rules
$$
t_0^2 = t_0,~~\tau_i^2 = 0,  ~~t_0 \tau_i = 0,~~\tau_i t_0 = \tau_i,~~i = 1, 2, \ldots.
$$
Hence, if $a \in I_\Bbbk(T)$, then $x a = 0$ and equation does not have solutions. If $a = \alpha t_0 + a_1$, where $\alpha \not = 0$, $a_1 \in I_\Bbbk(T)$, then 
$$
x a = x (\alpha t_0 + a_1) = \alpha x,
$$
and we have a solution $x = \alpha^{-1} b$.
\end{proof}

\bigskip

\section{Self-distributive (bi)algebras and related structures}
\label{self-rel}
\subsection{Linear racks and quandles and YBE}
Notions of the linear quandle and rack are very close to our rack and quandle bialgebras.
The following definition generalizes Definition 3.4  from \cite{XS}.

\begin{definition}
A counital coalgebra $(A; *, \Delta, \varepsilon)$ with a coalgebra morphism $* \colon A \otimes A \to A$ is said to be a {\it linear quandle} if  the following axioms hold:

1) Any element $a \in A$ is a generalized idempotent: 
 $$
\sum_i a^{(1)}_i * a^{(2)}_i = a.
$$

2) There exists a binary coalgebra morphism $\bar{*} \colon A \otimes A \to A$ which is right inverse to $*$: for any $a, b \in A$ it holds that
$$
(b * a^{(2)}) \, \bar{*} \, a^{(1)} = \varepsilon(a) b =(b \, \bar{*} \, a^{(2)}) *  a^{(1)}.
$$

3) Generalized self-distributivity:
$$
(a * b) * c = (a * c^{(1)}) * (b * c^{(2)}),
$$
for any $a, b, c \in A$.
\end{definition}

If $(A; *, \Delta, \varepsilon)$ is an algebra and a coassociative counital coalgebra which satisfies  axiom 3), then it is  called a linear shelf  \cite{XS}. If in additionally $(C; *, \Delta, \varepsilon)$ satisfies   axiom 2) , then it called  a linear rack \cite{XS}. We will call a linear rack satisfying axiom 1) a linear quandle. 

\begin{example}[\cite{ABRW}]
Let $(X, * )$ be a rack, $\Bbbk[X]$ be a rack algebra. This algebra equipped with the following coproduct 
$$
\Delta \colon \Bbbk[X] \to \Bbbk[X] \otimes \Bbbk[X],~~\Delta(x) = x \otimes x, ~~x \in X,
$$
and counit $\varepsilon \colon \Bbbk[X]  \to \Bbbk$, $\varepsilon (x) = 1$, $x \in X$, is a cocommutative coassociative counital coalgebra. Then $(\Bbbk[X]; *, \Delta, \varepsilon)$
 is a linear rack.
\end{example}

\begin{example}[\cite{L}] Let $(L; [\cdot, \cdot])$ be a Leibniz algebra, that is an algebra with an operation $[\cdot,\cdot]$ (non necessarily skew-symmetric) satisfying the Jacoby identity,
\bea
[[x,y],z]+[[y,z],x]+[[z,x],y]=0.\nn
\eea
Consider the vector space  $N=\Bbbk\oplus L$ and define two linear maps $\Delta \colon N \to N \otimes N$ by
\bea
\Delta(1)=1\otimes 1,\qquad \Delta(x)=x\otimes 1+1\otimes x,~~x \in L, \nn
\eea
and $\varepsilon \colon N \to \Bbbk$ by
$$
\varepsilon(1) = 1,~~\varepsilon(x)= 0,~~x \in L.
$$
It is evident that $(N; \Delta, \varepsilon)$ is a cocommutative coassociative counital coalgebra. Define the operation $*, \bar{*} \colon N \otimes N \to N$ as follows:
$$
1 * 1 = 1,~~1 * x = 0,~~x * 1 = x,~~x * y = [x, y];
$$
$$
1 \, \bar{*} \, 1 = 1,~~1 \, \bar{*} \, x = 0,~~x \, \bar{*} \, 1 = x,~~x \, \bar{*} \, y = - [x, y].
$$
Then   $(N; *, \Delta, \varepsilon)$ is a linear rack.
\end{example}

\begin{example} 
Let $(Q, * )$ be a quandle, $\Bbbk[Q]$ be a quandle algebra. This algebra equipped with the following coproduct 
$$
\Delta \colon \Bbbk[Q] \to \Bbbk[Q] \otimes \Bbbk[Q],~~\Delta(x) = x \otimes x, ~~x \in Q,
$$
and counit $\varepsilon \colon \Bbbk[Q]  \to \Bbbk$, $\varepsilon (x) = 1$, $x \in Q$, is a cocommutative coassociative counital coalgebra. Then $(\Bbbk[Q]; *, \Delta, \varepsilon)$
 is a linear quandle.
\end{example}

There is a connection of the linear racks with solutions of the Yang-Baxter equation:

\begin{proposition}[\cite{L}]
If $(A; *, \Delta, \varepsilon)$ is a cocommutative linear rack, then $R \colon A \otimes A \to A$ defined by
$$
R (a \otimes b) = b^{(1)} \otimes (a * b^{(2)}), ~~a, b \in A,
$$
satisfies the equation
$$
(R \otimes \id) (\id \otimes R) (R \otimes \id) = (\id \otimes R) (R \otimes \id) (\id \otimes R).
$$
Moreover, $R$ is invertible and its inverse is
$$
R^{-1} (a \otimes b) = (b \, \bar{*} \, a^{(1)}) \otimes a^{(2)}, ~~a, b \in A.
$$ 

\end{proposition}

The following example is Example 2.2 in \cite{ABRW}.

\begin{example}[\cite{ABRW}]
Let $(H, \Delta_H, \epsilon_H,  \mu_H, {\bf 1}_H, S)$ be a cocommutative Hopf algebra over $\Bbbk$, then the new multiplication
$$
\mu \colon H \otimes H \to H,~~\mu( h' \otimes h) = h' * h,
$$
defined by the usual adjoint representation
$$
h' * h = ad_h(h') = \sum h^{(1)}_i h' (S(h^{(2)}_i)),
$$
equips the $C^4$-coalgebra $(H, \Delta_H, \epsilon_H,  {\bf 1}_H)$ with a self-distributive bialgebra structure. In general, the adjoint representation does not seem to preserve the coalgebra structure if no cocommutativity is assumed.
\end{example}

\subsection{Quandle bialgebras and generalized Jordan algebras}
A Jordan algebra is a vector space $A$ with a commutative binary operation
\bea
m: A\otimes A\rightarrow A,\qquad m(a,b)=a\circ b,\nn
\eea
satisfying a Jordan property
\bea
\label{Jordan}
((a\circ a)\circ b)\circ a=(a\circ a)\circ(b\circ a).
\eea
This relation is fulfilled in quandles. Namely it is always true that 
\bea
a\circ a& =&a;\nn\\
(a\circ b)\circ a&=&a\circ(b\circ a).\nn
\eea
However the last properties of quandles are not relevant with the linear structure. We propose here a generalization for the Jordan algebra concept including a subclass of quandle algebras.

\begin{definition}
A generalized Jordan bialgebra is a vector space $A$ with compatible multiplication $m:A\otimes A\rightarrow A$ and  comultiplication $\Delta:A\rightarrow A\otimes A$ such that $\Delta$ is coassociative and both operations are related by
\bea
\label{gen-Jor}
((a^{(1)}\circ a^{(2)})\circ b)\circ a^{(3)}=(a^{(1)}\circ a^{(2)})\circ(b\circ a^{(3)}),
\eea
where we have used the Sweedler notations for the comultiplication
\bea
(\Delta\otimes \id)\Delta a=\sum_i a^{(1)}_i\otimes a^{(2)}_i\otimes a^{(3)}_i=a^{(1)}\otimes a^{(2)}\otimes a^{(3)}.\nn
\eea
\end{definition}
We get the following obvious
\begin{proposition}
Any quandle bialgebra is a generalized Jordan bialgebra.
\end{proposition}

\begin{proof} Let us see the compatibility condition for $(A; \circ, \Delta):$  
\bea
\Delta(a \circ b)=\Delta(a) \circ \Delta(b).
\eea
Let $$a=\sum^{n}_{i} \alpha_{i} q_{i} ,\,\,\,\,  b=\sum^{n}_{j} \beta _{j} q_{j},$$
where $\alpha_{i},  \beta _{j}, q_i, q_j \in \Bbbk$. Then the LHS can be expressed as
$$\Delta(a \circ b)=\Delta \left( \sum^{n}_{i,j} \alpha_{i}  \beta _{j} q_{i} q_{j}  \right)=  \sum^{n}_{i,j} \alpha_{i}  \beta _{j} \Delta \left( q_{i} q_{j}\right)=
 \sum^{n}_{i,j} \alpha_{i}  \beta _{j} \Delta \left( q_{i}q_{j} \otimes q_{i}q_{j}\right).$$
 The RHS is
 $$\Delta(a) \circ \Delta(b)=  \left( \sum^{n}_{i} \alpha_{i} q_{i} \otimes q_{i} \right) \circ  \left( \sum^{n}_{j}  \beta _{j} q_{j} \otimes q_{j} \right)=
  \sum^{n}_{i,j} \alpha_{i}  \beta _{j}  \left( q_{i} \otimes q_{i} \right)  \left( q_{j} \otimes q_{j} \right)= $$
 $$= \sum^{n}_{i,j} \alpha_{i}  \beta _{j} \Delta \left( q_{i}q_{j} \otimes q_{i}q_{j}\right).$$
Both sides of the equality \ref{gen-Jor} are bilinear on $a,b$ hence it could be verified on generators of the quandle bialgebra, but in this case is it just the equality \ref{Jordan}.
\end{proof}
\bigskip

\section{Classification of  2-dimension self-distributive bialgebras} \label{class}

\medskip

In this section we find  2-dimensional counital self-distributive bialgebras over the  field  $\Bbbk$ of $char (\Bbbk) \not=2$. In the case  of the field of complex numbers we give a full list of these bialgebras. 

Let $A$ be  a  vector space  over a field $\Bbbk$ with a basis $e_1$, $e_2$. We define  two operations: a multiplication $ A \otimes A \to A$ and a comultiplication 
$\Delta  \colon A \to  A \otimes A$ on $A$, such that the following axioms hold:
~\\
1) Coassociativity,
\bea
(\Delta\otimes \id_A)\Delta=(\id_A \otimes \Delta)\Delta.\nn
\eea
~\\
2) Generalized self-distributivity,
$$
(a b) c = (a c^{(1)})  (b c^{(2)}),~~~\Delta(c)=\sum_i c_i^{(1)} \otimes c_i^{(2)} =c^{(1)}\otimes c^{(2)}
$$
for any $a, b, c \in A$.
~\\
3) Compatibility  condition
\bea
\Delta(a b)=\Delta (a) \Delta(b)\nn
\eea
where the product on the right hand side is the component-vice product on $A\otimes A.$
\\
4) Counit, $\varepsilon \colon A \to \Bbbk$ such that
$$
(\varepsilon \otimes \id_A) \Delta =(\id_A \otimes \varepsilon) \Delta.
$$

\subsection{2-dimension counital coalgebras} 
We start by recalling the classification for the dual object $A^*$ - a 2-dimensional associative unital algebra. Over the field $\mathbb{C}$ the well known result says that if we introduce a basis $e_1=1$, $e_2$, then there are two non-isomorphic algebras with  multiplications 
\bea
1) && e_1 e_1 = e_1,~e_1 e_2 = e_2 e_1  = e_2,~e_2 e_2 = e_1;\nn\\
2) && e_1 e_1 = e_1,~e_1 e_2 = e_2 e_1 =  e_2,~e_2 e_2 = 0. \nn
\eea
\begin{remark}
The first algebra is equivalent to the algebra with multiplication
\bea
1')~ e'_1 e'_1 = e'_1,~e'_1 e'_2 = e'_2 e'_1  = 0,~e'_2 e'_2 = e'_2.\nn
\eea
Indeed, it is sufficient to put $e_1 = e_1' + e_2'$, $e_2 = e_1' - e_2'$.  
\end{remark}

\begin{remark}
In a case of a field $\Bbbk$ of $char (\Bbbk) \not=2$ a full description of 2-dimensional associative unital algebras  can be found in \cite[\S~1.5]{P}. In this case we have algebras 1), 2) and an algebra which is a quadratic extension of $\Bbbk$ by an element $e_2 \not\in \Bbbk$ which is a solution of a quadratic  equation $x^2 = a$. In this case we have a multiplication table
\bea
3)~ e_1 e_1 = e_1,~ e_1 e_2 = e_2 e_1  = e_2,~e_2 e_2 = a e_1.\nn
\eea
\end{remark}

\begin{example}
Let $\Bbbk = \mathbb{Q}$ be a field of rational numbers, $d \in \mathbb{Q}$ is square free, then $\mathbb{Q} \cdot 1 + \mathbb{Q} \cdot \sqrt{d}$ is a quadratic extension of $\mathbb{Q}$.
\end{example}

Let us consider the (dual) coassociative 2-dimension coalgebra $A = \langle f_1, f_2 \rangle$ where $f_i$ is the dual basis to $e_j$. 

\begin{lemma}
The comultiplication dual to the multiplication 1)  has the form
\bea
\Delta f_1 = f_1 \otimes f_1 + f_2 \otimes f_2,~~\Delta f_2 = f_1 \otimes f_2 + f_2 \otimes f_1.\nn
\eea
The comultiplication dual to the multiplication 2)  has the form
\bea
\Delta f_1 = f_1 \otimes f_1,~~\Delta f_2 = f_1 \otimes f_2 + f_2 \otimes f_1.\nn
\eea
The comultiplication dual to the multiplication 3)  has the form
\bea
\Delta f_1 = f_1 \otimes f_1 + a (f_2 \otimes f_2),~~\Delta f_2 = f_1 \otimes f_2 + f_2 \otimes f_1.\nn
\eea
\end{lemma}
\begin{proof}
Let us verify the statement for the type 1) comultiplication (proofs for comultiplications 2) and 3) are similar).
By the definition of  the dual basis we have
$$
f_1(e_1) = 1, f_1(e_2) = 0,~~f_2(e_1) = 0, f_2(e_2) = 1.
$$
Further, since
$$
\Delta f_k (e_i, e_j) = f_k (e_i e_j),~~~i, j, k \in \{ 1, 2 \},
$$
we have
$$
\Delta f_1 (e_1,  e_1) = f_1(e_1) = 1,~~\Delta f_1 (e_1,  e_2) = f_1(e_2) = 0,
$$
$$
\Delta f_1 (e_2,  e_1) = f_1(e_2) = 0,~~\Delta f_1 (e_2,  e_2) = f_1(e_1) = 1.
$$
The similar formulas for $f_2$ produce: 
$$
\Delta f_2 (e_1,  e_1) = f_2(e_1) = 0,~~\Delta f_2 (e_1,  e_2) = f_2(e_2) = 1,
$$
$$
\Delta f_2 (e_2,  e_1) = f_2(e_2) = 1,~~\Delta f_2 (e_2,  e_2) = f_2(e_1) = 0.
$$
\end{proof}

\begin{lemma}
\label{group_like}
The comultiplication $\Delta \colon A \to A \otimes A$, which is defined by the formulas
$$
\Delta f_1 = f_1 \otimes f_1 + f_2 \otimes f_2,~~\Delta f_2 = f_1 \otimes f_2 + f_2 \otimes f_1
$$
is equivalent to the group comultiplication:
$$
\Delta g_1 = g_1 \otimes g_1,~~\Delta g_2 =  g_2 \otimes g_2.
$$
\end{lemma}

\begin{proof}
Let us take other basis in $A$,
$$
g_1 = f_1 + f_2,~~g_2 = f_1 - f_2.
$$
Then
$$
\Delta g_1 = \Delta f_1 + \Delta f_2 = f_1 \otimes f_1 + f_2 \otimes f_2 +  f_1 \otimes f_2 + f_2 \otimes f_1 = (f_1 + f_2) \otimes (f_1 + f_2) = g_1 \otimes g_1,
$$
$$
\Delta g_2 = \Delta f_1 - \Delta f_2  = f_1 \otimes f_1 + f_2 \otimes f_2 -  f_1 \otimes f_2 - f_2 \otimes f_1 = (f_1 - f_2) \otimes (f_1 - f_2) = g_2 \otimes g_2.
$$
\end{proof}

\subsection{2-dimension self-distributive bialgebras  of type 1} 
\label{sec_type1}
Due to Lemma \ref{group_like} the comultiplication of the first type is equivalent to the following one:
\bea
\label{group-like-comult}
\Delta x = x \otimes x,~~\Delta y =  y \otimes y.
\eea
\begin{theorem}
Up to the permutation  of basis vectors each self-distributive 2-dimensional bialgebra with the comultiplication (\ref{group-like-comult}) can be represented as one of the list:
\bea
1) &&x^2 = x,\qquad x y = y, \qquad y x = x, \qquad y^2 = y;\nn\\
2) &&x^2 = x,\qquad x y = x, \qquad y x = y, \qquad y^2 = y;\nn\\
3) &&x^2 = x,\qquad x y = x, \qquad y x = x, \qquad y^2 = y;\nn\\
4) &&x^2 = x,\qquad x y = 0, \qquad y x = 0, \qquad y^2 = y;\nn\\
5) &&x^2 = x,\qquad x y = x, \qquad y x = y, \qquad y^2 = x;\nn\\
6) &&x^2 = x,\qquad x y = x, \qquad y x = x, \qquad y^2 = x;\nn\\
7) &&x^2 = x,\qquad x y = 0, \qquad y x = y, \qquad y^2 = 0;\nn\\
8) &&x^2 = x,\qquad x y = 0, \qquad y x = 0, \qquad y^2 = 0;\nn\\
9) &&x^2 = y,\qquad x y = y, \qquad y x = x, \qquad y^2 = x;\nn\\
10) &&x^2 = y,\qquad x y = y, \qquad y x = 0, \qquad y^2 = 0;\nn\\
11) &&x^2 = y,\qquad x y = 0, \qquad y x = 0, \qquad y^2 = 0;\nn\\
12) &&x^2 = 0,\qquad x y = y, \qquad y x = 0, \qquad y^2 = 0;\nn\\
13) &&x^2 = 0,\qquad x y = 0, \qquad y x = 0, \qquad y^2 = 0.\nn
\eea
\end{theorem}
\begin{proof}
The generalized self-distributivity condition in the case  of the group-like comultiplication is the following condition on the basis elements:
\begin{equation} \label{SD}
(a b) c = (a c^{(1)}) (b c^{(2)})= (a c)(b c).\nn
\end{equation} 
The compatibility  condition for the multiplication and comultiplication is the condition that the comultiplication $\Delta$ is an algebra homomorphism. This is equivalent to the following identity,
\bea
\Delta (a b) = \Delta(a) \Delta(b), ~~a, b \in A.\nn
\eea
Let us introduce a notation for the multiplication:
\bea
x x = a_1 x + a_2 y,~~x y = b_1 x + b_2 y,~~y x = c_1 x + c_2 y,~~yy = d_1 x + d_2 y.\nn
\eea
The compatibility  condition for $a=b=x$ gives
\bea
a_1 x \otimes x + a_2 y \otimes y = a_1^2 x \otimes x + a_1 a_2 x \otimes y + a_2 a_1 y \otimes x + a_2^2 y \otimes y.\nn
\eea
This immediately implies
\bea
a_1 a_2 = 0,\qquad a_1^2 = a_1,\qquad  a_2^2 = a_2.\nn
\eea
It means that the only possible pairs for $(a_1, a_2)$ are 
\bea
(a_1, a_2) \in \{ (0, 0),  (1, 0), (0, 1)\}.\nn
\eea
By the same arguments for $\Delta(xy),~\Delta(yx),~\Delta(y^2)$ we get the same conditions for other coefficients
$$
(d_1, d_2) \in \{ (0, 0),  (1, 0), (0, 1)\};
$$
$$
(b_1, b_2) \in \{ (0, 0),  (1, 0), (0, 1)\};
$$
$$
(c_1, c_2) \in \{ (0, 0),  (1, 0), (0, 1)\}.
$$
Then using the self-distributivity condition we get several properties. For example the equation
\bea
(x x)x=(x x)(x x)\nn
\eea
provides an implication
\bea
(a_1,a_2)=(0,1) \Rightarrow (d_1,d_2)=(c_1,c_2).\nn
\eea
The similar calculation for 
\bea
(y y)y=(y y)(y y)\nn
\eea
gives the analogous condition for  $d_i:$
\bea
(d_1,d_2)=(1,0) \Rightarrow (a_1,a_2)=(b_1,b_2).\nn
\eea
The generalized self-distributivity condition
\bea
(x x)y=(x y)(x y)\nn
\eea
produces the following
\bea
(a_1,a_2)=(1,0) \Rightarrow b_2 d_1=0;\quad b_2(d_2-1)=0.\nn
\eea
Considering all generalized self-distributivity conditions 
 and using the comultiplication symmetry $x\mapsto y; y\mapsto x$ we get the list of possible multiplications.
\end{proof}

\begin{remark}
This result corresponds to Lemma 3.8 from \cite[Section 3.4]{CCES}. 
\end{remark}

\bigskip

\subsection{2-dimension self-distributive bialgebras  of type 2} 

The second counital comultiplication is given by 
\bea
\label{nilp-comult}
\Delta x = x \otimes x,~~\Delta y =  x \otimes y+ y \otimes x.
\eea
From the condition $\Delta (x x) = \Delta(x) \Delta(x)$ it follows that
$$
a_1 x \otimes x + a_2 x \otimes y+a_2  y \otimes x = a_1^2 x \otimes x + a_1 a_2 x \otimes y + a_2 a_1 y \otimes x + a_2^2 y \otimes y.
$$
Hence, $a_1 a_2 = a_2$, $a_1^2 = a_1$, $a_2^2 = 0$, $a_1 a_2 = 0$.  $a_2 = 0,$ then $a_1 = 0$ or $a_1 = 1$.  We have two solutions:
$$
(a_1, a_2) \in \{ (0, 0),  (1, 0)\}.
$$
From the condition $\Delta (yy) = \Delta(y) \Delta(y)$ we get
$$
(d_1, d_2) \in \{ (0, 0),  (0, d_2)\}.
$$
From the condition $\Delta (x y) = \Delta(x) \Delta(y)$ we get
$$
(b_1, b_2) \in \{ (0, 0),  (0, b_2)\}.
$$
By the symmetry, from the condition $\Delta (y x) = \Delta(y) \Delta(x)$ we get
$$
(c_1, c_2) \in \{ (0, 0),  (0, c_2)\}.
$$

By the similar arguments as in Subsection \ref{sec_type1} we get:
\bigskip
\begin{theorem}
Each self-distributive 2-dimensional bialgebra with the comultiplication (\ref{nilp-comult}) can be represented as one of the list:
\bea
1) && x^2 = 0,\qquad x y = 0, \qquad y x = 0, \qquad y^2 = 0;\nn\\
2) && x^2 = 0,\qquad x y = y, \qquad y x = 0, \qquad y^2 = 0;\nn\\
3) && x^2 = x,\qquad x y = 0, \qquad y x = 0, \qquad y^2 = 0;\nn\\
4) && x^2 = x,\qquad x y = y, \qquad y x = 0, \qquad y^2 = 0.\nn\\
\eea
\end{theorem}

\subsection{2-dimension self-distributive bialgebras  of type 3}
Let us consider coalgebra with comutiplication
\bea
\label{comult_a}
\Delta x = x \otimes x + a (y \otimes y),~~\Delta y = x \otimes y + y \otimes x.
\eea

\begin{remark}
If $a = -1$, then the full list of such algebras can be found  in \cite[Section~3.4]{CCES}. 
\end{remark}

The condition $\Delta (x x) = \Delta(x) \Delta(x)$ implies
$$
a_1 x \otimes x + a_2 x \otimes y+a_2  y \otimes x + a a_1 y \otimes y = (a_1^2 + a c_1^2 + a b_1^2 + a^2 d_1^2) x \otimes x +
$$
$$
+ (a_1 a_2 + a c_1 c_2 + a b_1 b_2 + a^2 d_1 d_2 )x \otimes y + 
(a_1 a_2 + a c_1 c_2 + a b_1 b_2 + a^2 d_1 d_2 ) y \otimes x + 
(a_2^2 + a  c_2^2 + a  b_2^2 + a^2  d_2^2) y \otimes y.
$$
Hence, we have a system
$$
\begin{cases}
a_1 = a_1^2 + a c_1^2 + a b_1^2 + a^2 d_1^2,\\
a_2 =a_1 a_2 + a c_1 c_2 + a b_1 b_2 + a^2 d_1 d_2, \\
a a_1 = a_2^2 + a  c_2^2 + a  b_2^2 + a^2  d_2^2.
\end{cases}
$$

Moreover the equality $\Delta (y y) = \Delta(y) \Delta(y)$ produces
$$
d_1 x \otimes x + d_2 x \otimes y+d_2  y \otimes x + a d_1 y \otimes y = (2 a_1 d_1 + 2 b_1 c_1) x \otimes x +
$$
$$
+ (a_1 d_2 + b_1 c_2 + c_1 b_2 +  d_1 a_2 )x \otimes y + 
(a_2 d_1 + b_2 c_1 + b_1 c_2 +  a_1 d_2) y \otimes x + 
(2 a_2 d_2 + 2 b_2 c_2) y \otimes y.
$$
Hence, we get a system
$$
\begin{cases}
d_1 = 2 a_1 d_1 + 2 b_1 c_1,\\
d_2 = a_2 d_1 + b_2 c_1 + b_1 c_2 +  a_1 d_2, \\
a d_1 = 2a_{2} d_{2} +2b_{2}c_{2}.
\end{cases}
$$

From the condition $\Delta (x y) = \Delta(x) \Delta(y)$ it follows that
$$
b_1 x \otimes x +  b_2 x \otimes y+b_2  y \otimes x + a b_1 y \otimes y = (2 a_1 b_1 + 2 a c_1 d_1) x \otimes x +
$$
$$
+ (a_1 b_2 +  a_2 b_1 + a c_1 d_2 + a c_2 d_1 )x \otimes y + 
(a_1 b_2 +  a_2 b_1 + a c_1 d_2 + a c_2 d_1 ) y \otimes x + 
(2a_2 b_2 + c_2 d_2 + a c_2 d_2) y \otimes y.
$$
So, we have a system
$$
\begin{cases}
b_1 = 2 a_1 b_1 + 2 a c_1 d_1,\\
b_2 =a_1 b_2 +  a_2 b_1 + a c_1 d_2 + a c_2 d_1, \\
a b_1 = 2a_2 b_2 + c_2 d_2 + a c_2 d_2.
\end{cases}
$$

The condition $\Delta (y x) = \Delta(y) \Delta(x)$ implies
$$
c_1 x \otimes x +  c_2 x \otimes y+c_2  y \otimes x + a c_1 y \otimes y = (2 a_1 c_1 + 2 a b_1 d_1) x \otimes x +
$$
$$
+ (a_1 c_2 +  a_2 c_1 + a b_1 d_2 + a b_2 d_1 )x \otimes y + 
(a_1 c_2 +  a_2 c_1 + a b_1 d_2 + a b_2 d_1 ) y \otimes x + 
(2a_2 c_2 + b_2 d_2 + a c_b d_2) y \otimes y.
$$
Thus, we have a system
$$
\begin{cases}
c_1 = 2 a_1 c_1 + 2 a b_1 d_1,\\
c_2 =a_1 c_2 +  a_2 c_1 + a b_1 d_2 + a b_2 d_1, \\
a c_1 = 2a_2 c_2 + b_2 d_2 + a b_2 d_2.
\end{cases}
$$

Continuing such an analysis we arrive at the following partial list of multiplications:
\bea
1) && x^2=x- \frac {\sqrt a} 2 y,~ xy= - \frac 1 2 y,~yx=\frac 1 2 y,~ y^2= - \frac 1 {2 \sqrt a} y;\nn\\
2) && x^2=x+ \frac {\sqrt a} 2 y,~ xy= - \frac 1 2 y,~yx=\frac 1 2 y,~ y^2=  \frac 1 {2 \sqrt a} y;\nn\\ 
3) && x^2=x,~x y =0,~y x=-y,~y^2=0;\nn\\
4) && x^2=0,~ xy=0,~ yx=0,~ y^2=0;\nn\\
5) && x^2=x- \sqrt a y,~ xy= 0,~yx=0,~ y^2= 0;\nn\\
6) && x^2=x + \sqrt a y,~ xy= 0,~yx=0,~ y^2= 0.\nn
\eea

\begin{ack}
The work on section 1 was supported by the Russian Science Foundation project No. 20-71-10110 (https://rscf.ru/en/project/23-71-50012/) which funds the work of DT at the Yaroslavl State University. The work on section 4 was supported by the Ministry of Science and Higher Education of the Russian Federation (Agreement No. 075-02-2025-1636).
The work on parts 2, 3  was  supported by  the Ministry of Science and Higher Education of Russia (agreement No. 075-02-2025-1728/2).

We thank S. A. Panasenko  who suggested an idea of  proving
Proposition \ref{psd}.
\end{ack}
\medskip

\paragraph{Conflict of Interest:}  The authors declare that they have no
 conflicts of interest.


\begin{thebibliography}{HD}



\bibitem{Drinfeld} 
V. G. Drinfeld, \textit{On some unsolved problems in quantum group theory}, Quantum groups (Leningrad, 1990), 1--8, Lecture Notes in Math., 1510, Springer, Berlin,~1992.

\bibitem{Yang}
C. N. Yang, \textit{Some exact results for the many-body problem in one dimension with repulsive
delta-function interaction.} Phys. Rev. Lett., 19 (1967),   1312--1315.

\bibitem{Bax} 
R. J. Baxter, \textit{Partition function of the eight-vertex lattice model}. Ann. Physics, 70  (1972), 193--228.

\bibitem{Kas} 
C. Kassel, \textit{Quantum groups},  Springer 1995. Graduate texts in mathematics; vol.~155.

\bibitem{KT} C. Kassel, V. Turaev, \textit{Braid groups},  Springer 2008.  Graduate texts in mathematics;  vol.~247. 

\bibitem{FRT}
L.~D.~Faddeev, N.~Yu.~Reshetikhin, and L.~A.~Takhtajan, \textit{Quantization of Lie groups and Lie algebras}, Academic
Press, Inc., Boston, MA 1 (1988), 129--139.

\bibitem{J}
V. F. R. Jones, \textit{Baxterization}, Internat. J. Modern Phys. A 6 (1991), 2035--2043.

\bibitem{FR}
R. Fenn,  C. Rourke, \textit{Racks and links in codimension two}, J. Knot Theory Ramifications, 1 (1992), 343--406.

\bibitem{AG} 
N.~Andruskiewitsch, M.~Gra\~na, \textit{From racks to pointed Hopf algebras}, Adv. Math., 178 (2003), 177--243.


\bibitem{Kin}
M. K. Kinyon,  \textit{Leibniz algebras, Lie racks, and digroups},  J. Lie Theory, 17, no. 1 (2007),  99--114.

 \bibitem{CCES}
J. Carter, A. Crans, M. Elhamdadi, M. Saito,  \textit{Cohomology of the Categorical Self-Distributivity},  J. Homotopy Relat. Struct.,  3,  no. 1  (2008), 13--63.




\bibitem{ABRW}
C. Alexandre, M. Bordemann, S. Rivi\'ere, 
F.  Wagemann,  \textit{Structure theory of Rack-Bialgebras}, J Generalized Lie Theory Appl., 10, no. 1 (2016), DOI: 10.4172/1736-4337.1000244.

%\bibitem{Agu}
% M. Aguiar, \textit{Infinitesimal Hopf algebras}, Contemporary Mathematics 267 (2000) 1--29.


%\bibitem{ACM}
% A. Anquela, T. Cortes, and F. Montaner, \textit{Nonassociative coalgebras},  Comm. Alg., 22, no. 12 (1994), 4693--4716.

\bibitem{BPS} 
V. G. Bardakov, I. B. S. Passi, and Mahender Singh, \textit{Quandle rings},  J. Algebra Appl. 18 (2019), no. 8, 1950157, 23 pp, https://doi.org/10.1142/S0219498819501573.

\bibitem{EFT} 
M. Elhamdadi, N. Fernando and B. Tsvelikhovskiy, \textit{Ring theoretic aspects of quandles}, J. Algebra 526 (2019), 166--187.

\bibitem{ENS}
M. Elhamdadi, B. Nunez and M. Singh, \textit{Enhancements of link colorings via idempotents of quandle
rings}, J. Pure Appl. Algebra 227, no. 10 (2923),  ID 107400.


\bibitem{BPS-1}
V. G. Bardakov, I. B. S. Passi and M. Singh, \textit{Zero-divisors and idempotents in quandle rings}, Osaka J. Math., 59 (2022), 611--637.


\bibitem{ENSS}
M. Elhamdadi, B. Nunez, M. Singh, D. Swain, \textit{Idempotents, free products and quandle coverings}, Int. J. Math. 34, no. 3  (2023), Article ID 2350011, 27 p.













%\bibitem{EMSZ} M. Elhamdadi, A. Makhlouf, S. Silvestrov and E. Zappala, \textit{The derivation problem for quandle algebras}, Internat. J. Algebra Comput.,  32 (2022), 985--1007. 

%[12]



%\bibitem{Eisermann} M. Eisermann, \textit{Yang-Baxter deformations of quandles and racks}, Algebr. Geom. Topol., 5 (2005), 537--562.

%\bibitem{Q} D. Quillen,  \textit{Rational Homotopy Theory},  Ann. Math. 90 (1969), 205--295.


%\bibitem{P}
%R. S. Pierce,  \textit{Associative Algebras}, Springer - Verlag, New York Heidelberg Berlin,~1982.


%\bibitem{Z}
%	A. B. Zamolodchikov,	\textit{Tetrahedra equations and integrable systems in three dimensional space}, Zh. Eksp. Teor. Fiz. 79 (1980) 641664. [English translation: Soviet Phys. JETP 52 (1980) 325--326].

%\bibitem{Z-1}	A. B. Zamolodchikov, 	\textit{Tetrahedron equations and the relativistic S-matrix of straight-strings in 2+1 dimensions}, Commun. Math. Phys., 79  (1981), 489--505.

%\bibitem{STF}
%E. K. Sklyanin, L. A. Takhtadzhyan, and L. D. Faddeev, \textit{The quantum inverse problem method}. I,
%Teor. Mat. Fiz. 40 (1979), no. 2, 194--220; English transl., Theor. and Math. Phys. 40,  no.~2 (1979).

%\bibitem{TF}
 %L. A. Takhtadzhyan and L. D. Faddeev, \textit{The quantum method for the inverse problem and the XYZ
%Heisenberg model}, Uspekhi Mat. Nauk 34 (1979), no. 5, 13--63; English transl. in Russian Math. Surveys 34, no.~5 (1979). 



%\bibitem{Veselov} A. P. Veselov, {\em Integrable maps,} Russian Math. Surveys, 46:5 (1991), 1-51

%\bibitem{Ves2} A. P. Veselov, {\em Yang-Baxter maps and integrable dynamics}, Phys.Lett.A 314 (2003) 214

%\bibitem{BS} V. V. Bazhanov, S. M. Sergeev, \textit{Yang-Baxter maps, discrete integrable equations and quantum groups}, Nucl. Phys. B 926 (2018) 509--543.

\bibitem{Joyce} 
D. Joyce, \textit{A classifying invariant of knots, the knot quandle}, J. Pure Appl. Algebra, 23 (1982), 37--65.

\bibitem{Mat} 
S. Matveev, \textit{Distributive groupoids in knot theor}, Mat. Sb. (N.S.),  119 (161), no.~1 (9), 1982, 78--88 (in Russian).

%\bibitem{Markl}M. Markl,  \textit{Operads and PROPs}. Handbook of Algebra, 2008, 87-140. doi:10.1016/s1570-7954(07)05002-4. 





%\bibitem{Sol} A. Soloviev, \textit{Non-unitary set-theoretical solutions to the quantum Yang-Baxter equation}, Math. Res. Lett.,  7, no. 5-6 (2000), 577--596.

%\bibitem{LV} V. Lebed, A. Vendramin, \textit{Homology of left non-degenerate set-theoretic solutions to the Yang-Baxter equation}, Advances Math., 304 (2017), 1219--1261

%\bibitem{PT} M. M. Preobrazhenskaya, D. V. Talalaev, \textit{ Group extensions, fiber bundles, and a parametric Yang-Baxter equation}, Theoret. and Math. Phys., 207, no. 2 (2021), 670--677.


\bibitem{CP}
V. Chari and A. Pressley, \textit{A guide to quantum groups}, Cambridge University Press, Cambridge, 1995.


%\bibitem{BMP}

\bibitem{ZSSS}
K.~A.~Zhevlakov, A.~M.~Slinko, I.~P.~Shestakov, A.~I.~Shirshov, \textit{Rings that are nearly associative}, Academic Press, 1982.


\bibitem{XS}
N. Xu, Y. Sheng, \textit{The Yang--Baxter equation, Leibniz algebras, racks and relative algebraic structure}, arXiv:2410.15972.

\bibitem{L}
V. Lebed, \textit{Categorical aspects of virtuality and self-distributivity}, J. Knot Theory Ramifications, 22 (2013),
1350045.

\bibitem{P}
R. S. Pierce,  \textit{Associative Algebras}, Springer - Verlag, New York Heidelberg Berlin,~1982.


%\bibitem{RS} N. Y. Reshetikhin, M. A.  Semenov-Tian-Shansky,   \textit{ Quantum R-matrices and factorization problems.} Journal of Geometry and Physics, 5, no. 4  (1988), 533--550. doi:10.1016/0393-0440(88)90018-6 

%\bibitem{DrTw} V. Drinfeld,  \textit{Quasi-Hopf algebras},  Algebra i Analiz, 1, no. 6 (1989),  114--148.

%\bibitem{DrHam} V. G. Drinfeld,  \textit{Hamiltonian structures on Lie groups, Lie bialgebras and the geometric meaning of the classical Yang-Baxter equation}, Sov, Math, Dokl, 27 (1983), 68--71.

%\bibitem{Kul} P. Kulish, A. Mudrov,  \textit{On twisting solutions to the Yang-Baxter equation}, Czechoslovak Journal of Physics 50, no. 1 (2000), 115--122, DOI: 10.1023/A:1022885317520


%\bibitem{Leib} Sh. A. Ayupov; B. A. Omirov (1998). \textit{On Leibniz Algebras.} In Khakimdjanov, Y.; Goze, M.; Ayupov, Sh. (eds.). Algebra and Operator Theory Proceedings of the Colloquium in Tashkent, 1997. Dordrecht: Springer. pp. 1-13. ISBN 9789401150729.

% \bibitem{ER}  R.~Ehrenborg and M.~Readdy, \textit{Coproducts and the cd-index}, J. Algebraic Combin. 8 (1998) 273--299.

%\bibitem{JR} S.~A.~Joni and G.~C.`Rota, \textit{Coalgebras and Bialgebras in Combinatorics}, Studies in Applied Mathematics,  61 (1979), 93--139. Reprinted in Gian-Carlo Rota on Combinatorics: Introductory papers and commentaries (Joseph P.S. Kung, Ed.), Birkhauser, Boston (1995).

%\bibitem{Zot} A. Zotov

%\bibitem{G} M. E. Goncharov,  \textit{Structures of Malcev Bialgebras on a simple non-Lie Malcev algebra}, arXiv:1008.4214v1.


%\bibitem{Z} V.~N.~Zhelyabin,  \textit{Jordan bialgebras of symmetric elements and Lie bialgebras}, Siberian mathematical journal, 39, no. 2 (1998), 261-276.


%\bibitem{Z1} V.~N.~Zhelyabin,  \textit{Jordan bialgebras and their relation to Lie bialgebras}, Algebra and logic,  36, no. 1 (1997), 1--15.



%\bibitem{BN}
%V. Bardakov, T. Nasybullov,  \textit{Embeddings of quandles into groups, Journal of Algebra and Its Applications},  19, no. 7 (2020), 2050136 (20 pages).
%DOI: 10.1142/S0219498820501364.
%
%\bibitem{BarSin}
% V. Bardakov, M. Singh, \textit{Quandle cohomology, extensions and automorphisms}, J. Algebra, 585 (2021), 558--591.





%\bibitem{BMP}





\end{thebibliography}
\end{document}